\documentclass[12pt]{amsart}
\usepackage{latexsym,amsmath,amsfonts,amssymb,amsthm}
\def\Z{\mathbb Z}

\def\pmod #1{\ ({\rm mod}\ {#1})}
\newtheorem{theorem}{Theorem}[section]
\newtheorem{lemma}{Lemma}[section]
\newtheorem{corollary}{Corollary}[section]

\begin{document}
\title{On consecutive happy numbers}
\author{Hao Pan}
\address{Department of Mathematics, Nanjing University,
Nanjing 210093, People's Republic of China}
\email{haopan79@yahoo.com.cn}
\begin{abstract} Let $e\geqslant 1$ and $b\geqslant 2$ be integers.
For a positive integer $n=\sum_{j=0}^ka_j\times b^j$ with $0\leqslant a_j<b$, define
$$
T_{e,b}(n)=\sum_{j=0}^ka_j^e.
$$
$n$ is called $(e,b)$-happy if $T_{e,b}^r(n)=1$ for some
$r\geqslant 0$, where $T_{e,b}^r$ is the $r$-th iteration of
$T_{e,b}$. In this paper, we prove that there exist arbitrarily
long sequences of consecutive $(e,b)$-happy numbers provided that
$e-1$ is not divisible by $p-1$ for any prime divisor $p$ of
$b-1$.

\end{abstract} 
\maketitle

\section {Introduction}
\setcounter{equation}{0}
\setcounter{theorem}{0}
\setcounter{lemma}{0}
\setcounter{corollary}{0}

For an arbitrary positive integer $n$, let $T(n)$ be the sum of the squares of 10-adic digits of $n$. That is,
if we write
$n=\sum_{j=0}^ka_j\times10^j$ with $0\leqslant a_1,a_2,\ldots,a_k<10$, then
$T(n)=\sum_{j=0}^ka_j^2$. Also, let $T^r$ denote the $r$-th iteration of $T$, i.e.,
$$
T^r(n)=\underset{r \text{ times}}{\underbrace{T(T(\cdots T}}(n)\cdots)).
$$
In particular, we set $T^0(n)=n$. If $T$ is iteratively applied to $n$, it is easy to see (cf. \cite{Ho}) that we either get 1 or fall into
a cycle
$$
4\to 16\to 37\to 58\to 89\to 145\to 42\to 20\to 4.
$$
We say that $n$ is a happy number if we get 1 by applying $T$ to $n$ iteratively, i.e.,
$T^r(n)=1$ for some positive integer $r$.

Obviously here squares and 10-adic digits are very specialized. In general, we may replace the square by an arbitrary positive $e$-th power,
and replace the base 10 by an integer $b\geqslant2$. Let $T_{e,b}(n)$ be the sum of the $e$-th powers of $b$-adic digits of $n$, i.e.,
$$
T_{e,b}\bigg(\sum_{j=0}^ka_j\times b^j\bigg)=\sum_{j=0}^{k}a_j^e
$$
where $0\leqslant a_j<b$. And $n$ is called $(e,b)$-happy provided that there exists $r\geqslant 0$ such that
$T_{e,b}^r(n)=1$. Observe that
$$
T_{e,b}(n)<(b-1)^e(\log_b n+1).
$$
So if we iteratively apply $T_{e,b}$ to $n$, the process must reach some fixed points or cycles. And the fixed points and cycles of $T_{2,b}$ and $T_{3,b}$ ($2\leqslant b\leqslant 10$) have
been listed in \cite{GT}.

In the second edition of his famous book \cite[Problem E34]{Gu}, Guy asked whether there exists arbitrarily long sequences of consecutive (2,10)-happy numbers?
For example, the least five consecutive (2,10)-happy numbers are
$$
44488,\ 44489,\ 44490,\ 44491,\ 44492.
$$
In \cite{ES}, El-Sedy and Siksek gave an affirmative answer to Guy's question. As we will see soon, a key of El-Sedy and Siksek's proof is to find $h>0$ such that
$h+x$ is (2,10)-happy for each $x\in\{1, 4, 16, 20, 37, 42, 58, 89, 145\}$.

However, El-Sedy and Siksek's result cannot be extended to every power $e$ and base $b$. In fact, assume that $p$ is a prime divisor of $b-1$ and $e\equiv 1\pmod{p-1}$.
Then by Fermat's little theorem
$$
T_{e,b}\bigg(\sum_{j=0}^k a_j\times b^j\bigg)=\sum_{j=0}^k a_j^e\equiv\sum_{j=0}^k a_j\equiv\sum_{j=0}^k a_j\times b^j\pmod{p},
$$
that is, $T_{e,b}(n)\equiv n\pmod{p}$ for every $n$. Hence now $n$ is $(e,b)$-happy only if $n\equiv 1\pmod{p}$.
In this paper, we shall show that the above examples are the only exceptions.
\begin{theorem}
\label{t1}
Let $e\geqslant 1$ and $b\geqslant 2$ be integers. Suppose that $e\not\equiv 1\pmod{p-1}$ for any prime divisor $p$ of $b-1$. Then
for arbitrary positive integer $m$, there exists $l>0$ such that
$l+1,l+2,\ldots,l+m$ are all $(e,b)$-happy.
\end{theorem}

In view of Theorem \ref{t1}, we know that there exist arbitrarily
long sequences of consecutive $(2,b)$-happy numbers if $b$ is
even. For example, the least nine consecutive $(2,16)$-happy
numbers are
$$
65988605+i,\qquad i=0,1,\ldots,8.
$$
However, there exists no pair of consecutive $(2,12)$-happy
numbers less than $2^{32}-1$ (the maximal value of unsigned long
integers).

 The proof of Theorem \ref{t1} will be given in the next
section.

\section {Proof of Theorem \ref{t1}}
\setcounter{equation}{0}
\setcounter{theorem}{0}
\setcounter{lemma}{0}
\setcounter{corollary}{0}

Let $\Z^+$ denote the set of all positive integers. Since
$2-1=1\mid e-1$, below we always assume that $b$ is even. And for
convenience, we abbreviate `$(e,b)$-happy' to `happy' since $e$
and $b$ are always fixed. The following lemma is motivated by
El-Sedy and Siksek's proof in \cite{ES}.
\begin{lemma}
\label{l1}
Let $x$ and $m$ be an arbitrary non-negative integers.
Then for any $r\geqslant 1$, there exists a positive integer $l$
such that
$$
T_{e,b}^r(l+y)=T_{e,b}^r(l)+T_{e,b}^r(y)=x+T_{e,b}^r(y)
$$
for each $0\leqslant y\leqslant m$.
\end{lemma}
\begin{proof}
We make an induction on $r$. When $r=1$, choose a positive integer $s$ such that $b^s>m$ and let
$$
l_1=\sum_{j=0}^{x-1}b^{s+j}.
$$
Clearly
$$
T_{e,b}(l_1+y)=T_{e,b}(l_1)+T_{e,b}(y)=x+T_{e,b}(y)
$$
for any $0\leqslant y\leqslant m$.

Now assume $r>1$ and the assertion of Lemma \ref{l1} holds for the smaller values of $r$.
Since $T_{e,b}(n)\leqslant(b-1)^e(\log_bn+1)$,
there exists an $m'$ satisfying that
$T_{e,b}(y)\leqslant m'$ for all $0\leqslant y\leqslant m$.
Thus by the induction hypothesis, there exists an $l_{r-1}$ such that
$$
T_{e,b}^{r-1}(l_{r-1}+T(y))=T_{e,b}^{r-1}(l_{r-1})+T_{e,b}^{r-1}(T(y))=x+T_{e,b}^{r}(y).
$$
whenever $0\leqslant y\leqslant m$.

Let
$$
l_r=\sum_{j=0}^{l_{r-1}-1}b^{s+j}.
$$
where $s$ satisfies that $b^s>m$.
Then
$$
T_{e,b}^{r}(l_r)=T_{e,b}^{r-1}(T(l_{r}))=T_{e,b}^{r-1}(l_{r-1})=x,
$$
and for each $0\leqslant y\leqslant m$
\begin{align*}
T_{e,b}^{r}(l_r+y)=&T_{e,b}^{r-1}(T_{e,b}(l_r+y))=T_{e,b}^{r-1}(T_{e,b}(l_r)+T(y))\\
=&T_{e,b}^{r-1}(l_{r-1}+T(y))=T_{e,b}^{r-1}(l_{r-1})+T_{e,b}^{r}(y)=T_{e,b}^{r}(l_r)+T_{e,b}^{r}(y).
\end{align*}
\end{proof}

Suppose that a subset $D_{e,b}$ of positive integers satisfies that:

\medskip\noindent (1)\quad For any $n\in\Z^+$, there exists $r\geqslant 0$ such that $T_{e,b}^r(n)\in D_{e,b}$.

\medskip\noindent (2)\quad For any $x\in D_{e,b}$, $T(x)\in D_{e,b}$.

\medskip\noindent (3)\quad For any $x\in D_{e,b}$, there exists $r\geqslant 1$ such that $T_{e,b}^r(x)=x$.

\medskip\noindent
Then we say that $D_{e,b}$ is a cycle set for $T_{e,b}$. It is not difficult to see that $D_{e,b}$ is finite and uniquely determined by
$e$ and $b$. For example, $D_{2,10}=\{1, 4, 16, 20 37, 42, 58, 89, 145\}$.
\begin{corollary}
\label{c1}
Let $D_{e,b}$ is the cycle set for $T_{e,b}$. Assume that there exists $h\in\Z^+$ such that
$h+x$ is happy for any $x\in D_{e,b}$. Then for arbitrary $m\in\Z^+$, there exists $l\in\Z^+$ such that
$l+1,l+2,\ldots,l+m$ are all happy.
\end{corollary}
\begin{proof}
By the definition of cycle sets, there exists $r\in\Z^+$ such that $T_{e,b}^r(y)\in D$ for all $1\leqslant y\leqslant m$.
Applying Lemma \ref{l1}, we can find an $l\in\Z^+$ such that
$$
T_{e,b}^r(l+y)=h+T_{e,b}^r(y)
$$
whenever $1\leqslant y\leqslant m$. Thus by noting that $x$ is happy if and only if $T_{e,b}^r(x)$ is happy, we are done.
\end{proof}

However, in general, it is not easy to search such $h$ for $D_{e,b}$. With help of computers, when $e=2$ and $b=10$,
El-Sedy and Siksek found such
$$
h=\sum_{r=1}^{233192}9\times 10^{r+4}+20958
$$
by noting that $233192\times9^2+2^2+T_{2,10}(958+x)$ is happy for any $x\in D_{2,10}$. Fortunately, the following lemma
will reduce the requirement of $h$.
\begin{lemma}
\label{l2}
Let $D_{e,b}$ is the cycle set for $T_{e,b}$. Assume that for any $x\in D_{e,b}$, there exists $h_x\in\Z^+$ such that
both $h_x+1$ and $h_x+x$ are happy. Then there exists $h\in\Z^+$ such that
$h+x$ is happy for each $x\in D_{e,b}$.
\end{lemma}
\begin{proof}
We shall prove that under the assumptions of Lemma \ref{l2}, for any subset $S$ of $D_{e,b}$ with $1\in S$ and $|S|\geqslant 2$,
there exists $h_S\in\Z^+$ such that $h_S+x$ is happy for any $x\in S$.

The cases $|S|=2$ are trivial. Assume that $|S|>2$ and the assertion holds for any smaller value of $|S|$.
For any $1\not=x\in S$, since $h_x+1$ and $h_x+x$ are happy, there exists $r\in\Z^+$ such that
$$
T_{e,b}^r(h_x+1)=T_{e,b}^r(h_x+x)=1,
$$
and
$$
T_{e,b}^r(h_x+y)\in D_{e,b}\qquad\text{for all }y\in S
$$
by the definition of the cycle set. Let
$$
S^*:=\{T_{e,b}^r(h_x+y):\, y\in S\}.
$$
Then $1\in S^*\subseteq D_{e,b}$ and $|S^*|<|S|$. Thus by the induction hypothesis,
we can find $h_{S^*}\in\Z^+$ such that $h_{S^*}+T_{e,b}^r(h_x+y)$ is happy for any $y\in S$.
Also in view of Lemma \ref{l1}, there exists $l\in\Z^+$ satisfying that
$$
T_{e,b}^r(l+h_x+y)=h_{S^*}+T_{e,b}^r(h_x+y)
$$
provided that $y\in S$. It follows that $(l+h_x)+y$ is happy for any $y\in S$. All are done.
\end{proof}

\begin{lemma}
\label{l3}
Suppose that for any integer $a$, there exists a happy
number $h$ such that
$$
h\equiv a\pmod{(b-1)^e}.
$$
Then for any $x\in\Z^+$, there exists an arbitrarily large happy number $l$ such that
$l+x$ is also happy.
\end{lemma}
\begin{proof} Choose $s\in\Z^+$ satisfying that
$b^s>x$ and let $x^*=b^s-x$. Suppose that $h$ is the happy number such that
$$
h\equiv T_{e,b}(x^*)\pmod{(b-1)^e}.
$$
Note that $hb^{\phi((b-1)^e)}$ is also happy and
$$
hb^{\phi((b-1)^e)}\equiv h\pmod{(b-1)^e}
$$
where $\phi$ is the Euler totient function.
We may assume that $h>T_{e,b}(x^*)$. Write $h=(b-1)^ek+T_{e,b}(x^*)$. Let
$$
l=x^*+\sum_{j=0}^{k-1}(b-1)b^{s+j}.
$$
Then
$$
T_{e,b}(l)=k(b-1)^e+T_{e,b}(x^*)=h
$$
and
$$
T_{e,b}(l+x)=T_{e,b}(b^s+\sum_{j=0}^{k-1}(b-1)b^{s+j})=T_{e,b}(b^{s+k})=1.
$$
It follows that both $l$ and $l+x$ are happy.
\end{proof}
\begin{lemma}
\label{l4} Let $n$ be a positive odd integer. Then for any $a$
with $a\equiv 1\pmod{n}$ and positive integer $k$, there exists
$r\in\Z^+$ such that
$$
(n+1)^r\equiv a\pmod{n^k}.
$$
\end{lemma}
\begin{proof}
Assume that $n=p_1^{\alpha_1}p_2^{\alpha_2}\cdots p_s^{\alpha_s}$
where $p_1, p_2,\ldots p_s$ are distinct odd primes and $\alpha_1, \alpha_1, \ldots, \alpha_s\geqslant 1$.
For $1\leqslant i\leqslant s$, let $g_i$ be a prime root of
$p_i^{\alpha_ik}$. Assume that
$$
n+1\equiv g_i^{\beta_i}\pmod{p_i^{\alpha_ik}}\text{ and }a\equiv
g_i^{\gamma_i}\pmod{p_i^{\alpha_ik}}
$$
for each $1\leqslant i\leqslant s$. Clearly both $\beta_i$ and
$\gamma_i$ are divisible by $\phi(p_i^{\alpha_i})$ since $n+1\equiv a\equiv\pmod{p_i^{\alpha_i}} $. So we only
need to find $r$ satisfying that
$$
\beta_i r\equiv\gamma_i\pmod{\phi(p^{\alpha_ik})}
$$
for all $i$, or equivalently,
$$
(\beta_i/\phi(p_i^{\alpha_i}))r\equiv\gamma_i/\phi(p_i^{\alpha_i})\pmod{p^{\alpha_i(k-1)}}.
$$
Note that $p_i\nmid \beta_i/\phi(p_i^{\alpha_i})$ since
$n+1\not\equiv1\pmod{p_i^{\alpha_i+1}}$. Thus such $r$ always
exists in view of the Chinese remainder theorem.
\end{proof}
\begin{corollary}
\label{c2}
Assume that for any integer $a$, there exists a happy number $h$
such that
$$
h\equiv a\pmod{b-1}.
$$
Then we can find a happy number $h'$ such that
$$
h'\equiv a\pmod{(b-1)^e}.
$$
\end{corollary}
\begin{proof} Suppose
that
$$
\sum_{j=1}^{h-1}b^{j}\equiv k_1(b-1)+a-1\pmod{(b-1)^e}.
$$
And suppose that
$$
b^{-h}(-k_1(b-1)+1)\equiv k_2(b-1)+1\pmod{(b-1)^e}.
$$
In light of Lemma \ref{l4}, there exists $r\in\Z^+$ such that
$$
b^r\equiv k_2(b-1)+1\pmod{(b-1)^e}.
$$
Therefore
$$
\sum_{j=1}^{h-1}b^{j}+b^{h+r}\equiv\sum_{j=1}^{h-1}b^{j}+b^h(k_2(b-1)+1)\equiv
a\pmod{(b-1)^e},
$$
which is apparently happy.
\end{proof}
\begin{lemma}
\label{l5}  Let $a$ be a positive integer. Assume that there
exists a happy number $h$ such that
$$
l\equiv T_{e,b}(h')\pmod{b-1}
$$
for some $h'\equiv a\pmod{b-1}$. Then we can find a happy number
$h$ such that
$$
h\equiv a\pmod{b-1}.
$$
\end{lemma}
\begin{proof}
In view of the proof of Lemma \ref{l2}, we may assume that $l>T_{e,b}(h')$. And let
$$
h=\sum_{j=0}^{l-T_{e,b}(h')}b^{s+j}+h',
$$
where we choose $s$ such that $b^s>h'$. Clearly
$$
h\equiv l-T_{e,b}(h')+h'\equiv a\pmod{b-1}.
$$
And
$$
T_{e,b}(h)=l-T_{e,b}(h')+T_{e,b}(h')=l.
$$
Thus $h$ is the desired happy number.
\end{proof}

Note that the property of $e$ is not used until now.
\begin{proof}[Proof of Theorem \ref{t1}]
Write $b-1=p_1^{\alpha_1}p_2^{\alpha_2}\cdots p_s^{\alpha_s}$
where $p_1, p_2,\ldots, p_s$ are distinct odd primes and $\alpha_1,
\alpha_2,\ldots, \alpha_s$ are positive integers. Let $1\leqslant g_i\leqslant p_i^{\alpha_i}$ be a primitive root of $p_i^{\alpha_i}$ for $1\leqslant i\leqslant s$.
For every positive integer $0\leqslant a\leqslant b-1$, let $L(a)\in\{0,1,\ldots,b-1\}$ be the integer such that
$$
L(a)\equiv\begin{cases}
a-g_i+g_i^e\pmod{p_i^{\alpha_i}}&\text{if }a\not\equiv 1\pmod{p_i^{\alpha_i}},\\
1\pmod{p_i^{\alpha_i}}&\text{if }a\equiv 1\pmod{p_i^{\alpha_i}},
\end{cases}
\quad i=1, 2, \ldots, s.
$$
And let $r_a\geqslant 0$ be the minimal integer such that $L^{r_a}(a)\equiv 1\pmod{b-1}$, where $L^r$ denotes $r$-th iteration of $L$.
Since $e\not\equiv 1\pmod{p_i-1}$, we have $g_i^e-g_i$ is prime to $p_i$ for every $i$. Hence $r_a$ always exists.

Combining Corollary \ref{c1}, Lemma \ref{l2}, Lemma \ref{l3} and Corollary \ref{c2}, now it suffices to show that for every integer $0\leqslant a\leqslant b-1$ there exists a happy number $h$ such that
$$
h\equiv a\pmod{b-1}.
$$
We use an induction on $r_a$. If $r_a=0$, then $a\equiv 1\pmod{b-1}$. There is noting to do. Now assume that $r_a\geqslant 1$ and the assertion holds for any $a'$ with $r_{a'}<r_a$.
Clearly $r_{L(a)}=r_a-1$. Hence by the induction hypothesis, there exists a happy number $l$ such that
$$
l\equiv L(a)\pmod{b-1}.
$$
Let $0\leqslant g\leqslant b-1$ be the integer such that
$$
g\equiv\begin{cases}
g_i\pmod{p_i^{\alpha_i}}&\text{if }a\not\equiv 1\pmod{p_i^{\alpha_i}},\\
1\pmod{p_i^{\alpha_i}}&\text{if }a\equiv 1\pmod{p_i^{\alpha_i}},
\end{cases}
\quad i=1, 2, \ldots, s.
$$
And let
$$
h'=\sum_{j=1}^{a+b-1-g}b^j+g.
$$
Then
$$
h'\equiv a+b-1-g+g\equiv a\pmod{b-1}.
$$
And
$$
T_{e,b}(h')\equiv a+b-1-g+g^e\equiv\begin{cases}
a-g_i+g_i^e\pmod{p_i^{\alpha_i}}&\text{if }a\not\equiv 1\pmod{p_i^{\alpha_i}},\\
1\pmod{p_i^{\alpha_i}}&\text{if }a\equiv 1\pmod{p_i^{\alpha_i}},
\end{cases}
$$
for every $1\leqslant i\leqslant s$. Hence
$$
T_{e,b}(h')\equiv L(a)\equiv l\pmod{b-1}.
$$
Thus in light of Lemma \ref{l5}, we are done.
\end{proof}

\end{document}